\documentclass[12pt]{amsart}

\usepackage{amssymb}

\newcommand{\ol}{\overline}

\newcommand{\al}{\alpha}

\newcommand{\dt}{\delta}

\newcommand{\gm}{\gamma}

\newcommand{\Dt}{\Delta}

\newcommand{\Tht}{\Theta}

\newcommand{\ord}{\mathbb{\rm ord}}

\newcommand{\tm}{\times}

\newcommand{\wdt}{\widetilde}

\newcommand{\bC}{\mathbb{C}}

\newcommand{\calP}{\mathcal{P}}
\newcommand{\calA}{\mathcal{A}}
\newcommand{\calW}{\mathcal{W}}

\newcommand{\mbA}{\mathbf{A}}
\newcommand{\mba}{\mathbf{a}}
\newcommand{\mbB}{\mathbf{B}}
\newcommand{\mbb}{\mathbf{b}}
\newcommand{\mbu}{\mathbf{u}}
\newcommand{\mbv}{\mathbf{v}}

\newcommand{\mbV}{\mathbf{V}}

\begin{document}

\begin{center}
Numerical comparison of different algorithms for construction of
wavelet matrices
\\[5mm]
N. Salia, A. Gamkrelidze, and L. Ephremidze
\end{center}

\vskip+0.5cm

\noindent {\small {\bf Abstract.} Factorization of compact wavelet
matrices into primitive ones has been known for more than 20
years. This method makes it possible to generate wavelet matrix
coefficients and also to specify them by their first row.
Recently, a new parametrization of compact wavelet matrices of the
same order and degree has been introduced by the last author. This
method also enables us to fulfill the above mentioned tasks of
matrix constructions. In the present paper, we briefly describe
the corresponding algorithms based on two different methods, and
numerically compare their performance.}

 \vskip+0.2cm\noindent  {\small {\bf Keywords:} Wavelet matrices,
paraunitary matrix polynomials, wavelet matrix completion
algorithm.}

\vskip+0.2cm \noindent  {\small {\bf  AMS subject classification
(2010):} 42C40}

\vskip+0.5cm

\section{Introduction}

An $m\tm(N+1)m$ matrix
\begin{equation}
\label{wm} {\mathcal A}=(A_0\,A_1\,\ldots\,A_N)=
\left(\begin{matrix}
a^1_{1}&a^1_{2}&\cdots&a^1_{(N+1)m}\\[1mm]
a^2_{1}&a^2_{2}&\cdots&a^2_{(N+1)m}\\
\vdots&\vdots&\cdots&\vdots\\
a^{m}_{1}&a^{m}_{2}&\cdots&a^m_{(N+1)m}\\
\end{matrix}\right)
\end{equation}
($A_j$ are square blocks) is called a {\em wavelet matrix}
\cite{KT} if it satisfies the so called {\em shifted orthogonality
condition}:
\begin{equation}
\label{soc}\sum_{j=0}^{N-k}A_jA_{j+k}^*=\dt_{k0}I_m\,,\;\;\;\;k=0,1,\ldots,N,
\end{equation}
where $A^*$ denotes the conjugate transpose of $A$, $\dt_{k0}$ is
the Kronecker delta, and $I_m$ is the $m\tm m$ unit matrix.

In the {\em polyphase representation} of matrix $\calA$,
\begin{equation}
\label{pr} {\mathbf A}(z)=\sum_{k=0}^N
A_kz^k=:\{\mba_{ij}(z)\}_{i,j=1}^m\,,
\end{equation}
 the condition (\ref{soc})
is equivalent to
\begin{equation}
\label{fac} {\mathbf A}(z)\wdt{\mbA}(z)=I_m\,,
\end{equation}
where $\wdt{\mbA}(z)=\sum_{k=0}^N A_k^*z^{-k}$ is the {\em
adjoint} to $\mbA(z)$.

In the sequel, the matrices of the form (\ref{wm}) and their
polyphase representation (\ref{pr}) will be identified.

Our notion of a wavelet matrix is weaker than usual as some linear
condition is also required to be satisfied (see, e.g. \cite{RW})
which is irrelevant in our consideration. Instead, we require the
condition
\begin{equation}
\label{erti} {\mathbf A}(1)=I_m\,.
\end{equation}

The integers $m$ and $N$ are called, respectively, the {\em rank}
and the {\em order} of a wavelet matrix (\ref{wm}) or (\ref{pr})
(it is assumed that $A_N\not=\mathbf{0}$). It follows from
(\ref{fac}) that $\det\mbA(z)$ has always the form $cz^d$, $d\geq
0$, $|c|=1$, and the integer $d$ is called the {\em degree} of
$\calA$. The class of wavelet matrices of rank $m$, order $N$ and
degree $d$ will be denoted by $\calW(m,N,d)$. In addition,
$\calW_0(m,N,d)$ denotes the class of those $\calA\in\calW(m,N,d)$
for which (\ref{erti}) holds, and $\calW_1(m,N,d)$ denotes the
class of those $\calA\in\calW_0(m,N,d)$ for which the last row of
$A_N$ differs from the zero vector of $\bC^m$.

It can be proved that the degree of any wavelet matrix is grater
than or equal to its order, i.e. $d\geq N$ (see, e.g., \cite[Lemma
1]{EL}) and $d=N$ holds except for some degenerated cases (see
\cite[p. 58]{RW}). We call the case $d>N$ singular as the
uniqueness of solutions, which we are going to construct
numerically, fails to hold in this situation \cite{KT}. Namely, we
consider wavelet matrices from the class $\calW(m,N,N)$. It
differs from $\calW_1(m,N,N)$ by unitary multipliers on the left
and on the right (see \cite{EL}).

A wavelet matrix $\mbV(z)$ of  order and degree 1 is called
primitive. It can be shown (see, e.g. \cite[p. 59]{RW}, \cite{KT},
\cite{EL}) that every $\mbV(z)\in \calW_0(m,1,1)$ has the form
$$
\mbV(z)=I_m-\mbv^*\mbv-\mbv^*\mbv z,
$$
where $\mbv=(v_1,v_2,\ldots,v_m)\in \bC^m$ is a vector of the unit
norm, $\mbv\mbv^*=1$.

The following wavelet matrix factorization theorem was first
proved in the yearly 90's in a related theory of multirate filter
banks \cite{Vai}. We formulate it for nonsingular matrices

\smallskip

{\bf Theorem 1.1.} For any $\mbA(z)\in \calW_0(m,N,N)$, there
exists a unique factorization
\begin{equation*}
\mbA(z)=\prod_{j=1}^N\mbV_j(z),
\end{equation*}
where each $\mbV_j(z)\in \calW_0(m,1,1)$.

\smallskip

This theorem provides a possibility to generate wavelet matrices
of arbitrary order. The computational complexity of this method
and its numerical tests are described in the next sections.
Theorem 1.1. helps also to solve the following wavelet matrix
completion problem  \cite{Hel}, \cite{KT}: Given the first row of
a wavelet matrix, find its remaining rows, i.e. if the first row
of (\ref{wm}) is given which satisfies the shifted orthogonality
condition
\begin{equation}
\label{pstr}
\sum_{j=1}^{(N+1-k)m}a^1_j\ol{a}^1_{j+km}=\dt_{k0}\,,\;\;\;k=0,1,\ldots,N,
\end{equation}
then one should find the remaining entries of $\calA$ which
results in wavelet matrix. We emphasize that this problem has a
unique solution (up to certain unitary matrix) if we search
$\calA$ in $\calW(m,N,N)$ (see \cite[Th. 4.17]{RW}). In the next
sections, we describe and test numerically the existing algorithm
of such construction.

A new parametrization of nonsingular compact wavelet matrices
appeared in \cite{EL} in the form of Theorem 1.2 below which gives
a one-to-one continuous map between $\bC^{N(m-1)}$ and
$\calW_1(m,N,N)$

Let   $\calP_N^+:=\big\{\!\sum_{k=0}^N c_k z^k:c_0,c_1,\ldots,c_N$
$\in\bC\big\}$ be the set of polynomials and
$\calP_N^-:=\big\{\sum_{k=1}^N c_k
z^{-k}:c_1,c_2,\ldots,c_N\in\bC\big\}$  (note that
$\calP_N^+\cap\calP_N^-=\{0\}$ according to our notation). If
$p(z)=\sum_{k=-N}^N c_k z^k$, then $[p(z)]^-=\sum_{k=-N}^{-1} c_k
z^k$ and $[p(z)]^+=\sum_{k=0}^{N} c_k z^k$.

\smallskip

{\bf Theorem 1.2.} Let $N\geq 1$. For any polynomials
\begin{equation}
\label{zeta} \zeta_j(z)\in\calP_N^-, \;j=1,2,\ldots, m-1,
\end{equation}
there exists a unique
\begin{equation}
\label{Az1} A(z)\in\calW_1(m,N,N)
\end{equation}
such that
\begin{equation}
\label{P+}
\zeta_1(z)a_{1j}(z)+\zeta_2(z)a_{2j}(z)+\ldots+\zeta_{m-1}(z)a_{m-1,j}(z)
+\wdt{a_{mj}}(z)\in \calP_N^+,\;\;j=1,2,\ldots,m.
\end{equation}
Conversely, for each $A(z)$ satisfying $(\ref{Az1})$, there exists
a unique $(m-1)$-tuple of Laurent polynomials $(\ref{zeta})$ such
that $(\ref{P+})$ holds.

\smallskip

Further refinement of Theorem 1.2 enables us to solve the wavelet
matrix completion problem as well \cite[\S 5]{EL}. The exact
formulas of these constructions and numerical tests of
corresponding algorithms are given in Sections 3 and 4.

In conclusion we analyze numerical performances of described
algorithms and, based on these data, compare two different
methods.

\section{The existing algorithms of wavelet matrix construction}

The following wavelet matrix generation procedure is based on
Theorem 1.1

\smallskip

{\bf Algorithm 2.1.} {\em Step 1.} Take arbitrary nonzero vectors
$\mbv_j\in\bC^m$, $j=1,2,\ldots,N$, (they can be selected
randomly) and let
\begin{equation}
\label{Pj} P_j=(\mbv_j\mbv_j^*)^{-1}\mbv_j^*\mbv_j\,.
\end{equation}
Then $\mbV_j(z)=I_m-P_j+P_jz$, $j=1,2,\ldots,N$, are primitive
wavelet matrices.

{\em Step 2.} Let $\mbA_0(z)=I_m$ and for $j=1,2,\ldots,N$ do:
\begin{equation}
\label{st2} \mbA_j(z)=\mbA_{j-1}(z)\big(I_m-P_j(I_m-z)\big).
\end{equation}
(Matrix multiplication in (\ref{st2}) requires approximately
$m^2(j-1)$ operations (ops) counting only multiplications. Thus
the cycle in Step 2 needs $\approx\sum_{j=1}^Nm^2(j-1)=O(m^2N^2)$
ops.)

Then $\mbA(z)=\mbA_N(z)$ will be the wavelet matrix of rank $m$
and degree $N$. It can be seen that $\ord(\calA)=N$ if and only if
$\mbv_j\mbv^*_{j+1}\not=\mathbf{0}$ for $j=1,2,\ldots,N-1$ (i.e.
the consecutive $\mbv_j$-s in  (\ref{Pj}) are not orthogonal), and
in this case the last row of $A_N=\prod_{j=1}^NP_j$ differs from
$\mathbf{0}\in\bC^m$ if and only if the last coordinate of
$\mbv_1$ differs from $0$. Thus, for randomly selected $\mbv_j$-s
in Step 1, the wavelet matrix $\mbA(z)$ belongs to
$\calW_1(m,N,N)$ with probability 1.

\smallskip

The following procedures describe a numerical solution to the
wavelet matrix completion problem \cite{KT}, \cite{Hel}.

\smallskip

{\bf Algorithm 2.2.} Given
\begin{equation}
\label{a1} \mba=(a^1_{1}, a^1_{2}, \cdots,
a^1_{(N+1)m})=:(\mba_0,\mba_1\ldots,\mba_{N}),\;\;\;\mba_N\not=\mathbf{0}
\end{equation}
satisfying conditions (\ref{pstr}) and
$\sum_{i=0}^N\mba_i=e_1=(1,0,\ldots,0)\in\bC^m$.

{\em Step 1.} Let $P_N=(\mba_N\mba_N^*)^{-1}\mba_N^*\mba_N$ and
let
$(\mba^{(N)}_0,\mba^{(N)}_1\ldots,\mba^{(N)}_{N}):=(\mba_0,\mba_1\ldots,\mba_{N})$.
For $j=N, N-1,\ldots,2$ do:
$$
\mba^{(j-1)}_{i}=\mba^{(j)}_{i}+(\mba^{(j)}_{i+1}-\mba^{(j)}_{i})P_j,\;\;i=0,1,\ldots,N-1,
$$
and
$$
P_{j-1}=\big(\mba^{(j-1)}_{N-1}(\mba^{(j-1)}_{N-1})^*\big)^{-1}(\mba^{(j-1)}_{N-1})^*\mba^{(j-1)}_{N-1}\,.
$$
(This step needs approximately $O(mN^2)$ ops.)

{\em Step 2.} Compute the product
\begin{equation*}
 \mbA(z)=\prod_{j=1}^N(I_m-P_j+P_jz)
\end{equation*}
using Step 2 of Algorithm 2.1.

Then $\calA=\mbA(z)$ is the unique wavelet matrix from
$\calW_0(m,N,N)$ with the first row (\ref{a1}) (see also \cite[Th.
4.17]{RW}).

All in all, the number of operations in Algorithms 2.1 and 2.2 can
be estimated as $O(m^2N^2)$.

\section{New algorithms of wavelet matrix construction}

In this section we describe algorithms based on recently developed
method of wavelet matrix parametrization \cite{JLE}. First we
generate $\calA\in\calW_1(m,N,N)$ (see \cite{EL} for justification
of the given procedures).

\smallskip

{\bf Algorithm 3.1.} {\em Step 1.} Take arbitrary $m-1$ Laurent
polynomials from $\calP_N^-$
\begin{equation}
\label{zeta1} \zeta_i(z)=\sum_{k=1}^N \gm_{ik}z^{-k},\;\;\;
i=1,2,\ldots,m-1,
\end{equation}
(the coefficients $\gm_{ik}$ can be selected randomly).

{\em Step 2.} Perform upper triangular, diagonal, lower triangular
factorization
\begin{equation}
\label{udu} \Dt=UDU^*
\end{equation}
of
\begin{equation}
\label{Dt1} \Dt=\sum_{i=1}^{m-1}\Tht_i\ol{\Tht_i}+I_{N+1},
\end{equation}
where $\Tht_i$ is the upper triangular $(N+1)\tm(N+1)$ Hankel
matrix with the first row $(0,\gm_{i1},\gm_{i2},\ldots,\gm_{iN})$.

Since $\Dt$ has a displacement structure of rank $m$ (see
\cite[Appendix]{JLE}) the factorization (\ref{udu}) can be
performed in $O(mN^2)$ ops  (as it is described in \cite[Appendix
F.1]{Kai} ) without constructing  (\ref{Dt1}) explicitly

{\em Step 3.} Solve the system of $(N+1)\tm(N+1)$ linear algebraic
equations
\begin{equation}
\label{ls}\Dt X=B_j
\end{equation}
$m$ times taking different right hand sides, where
$B_j=(0,\gm_{i1},\gm_{i2},\ldots,\gm_{iN})^T$, $j=1,2,\ldots,m-1$,
and $B_m=(1,0,\ldots,0)^T$.

Since we have the factorization (\ref{udu}), the solution of the
system (\ref{ls}) requires  $O(N^2)$ ops and Step 3 totally needs
$O(mN^2)$ ops.

Let $(\al_{j0},\al_{j1},\ldots,\al_{jN})$ be the solution of
(\ref{ls}) and let
$$
\mbu_j(z)=\sum_{k=0}^N\al_{jk}z^{-k} \text{ and }
\mbb_{mj}(z)=z^N\mbu_j(z),\;\;\;\;j=1,2,\ldots,m.
$$

{\em Step 4.} Compute the coefficients of the following
polynomials from $\calP_N^+$
$$
\mbb_{ij}(z)=[\wdt{\zeta}_i(z)\mbu_j(z)]^+-\dt_{ij},\;\;\;1\leq
i<m,\;1\leq j\leq m.
$$
As the multiplication of polynomials of order $N$ takes $O(N\log
N)$ ops by FFT, Step 4 totally needs $O(m^2N\log N)$ ops.

{\em Step 5.} Constructing the matrix polynomial
$\mbB(z)=\{\mbb_{ij}(z)\}_{i,j=1}^m$,
$$
\mbA(z)=\mbB(z)\big(\mbB(1)\big)^{-1}
$$
will be a wavelet matrix from $\calW_1(m,N,N)$.

Since $m\tm m$ matrix inversion needs $O(m^3)$ ops, Step 5 totally
needs $O(Nm^3)$ ops.

\smallskip

Now we describe a new algorithm of wavelet matrix completion based
on Theorem 1.2. Its justification can be found in \cite{EL}.

\smallskip

{\bf Algorithm 3.2.} Data is the same as in Algorithm 2.2.

{\em Step 1.} Select a coordinate of $\mba_N=(a^1_{mN+1},
a^1_{mN+2}, \cdots, a^1_{m(N+1)})$  with maximum absolute value.
Since $\mba_N\not=\mathbf{0}$, this coordinate differs from $0$
and let it be $a^1_{mN+j}$. This preliminary step will improve the
accuracy of the final result.

{\em Step 2.} Let
$\big(\mba_{11}(z),\mba_{12}(z),\ldots,\mba_{1m}(z)\big)$ be the
polyphase representation of (\ref{a1}), i.e. the first row of
(\ref{pr}).

Compute the first $N+1$ coefficients, say
$\gm_0,\gm_1,\ldots,\gm_N$, of the reciprocal (in a neighborhood
of $0$) of
$\sum_{k=0}^N\ol{a}_{mk+j}^1z^{N-k}=z^N\wdt{\mba_{1j}}(z)$, where
$j$ was determined in Step 1. (This step requires $O(N^2)$ ops,
though some papers \cite{bini} report that it can be done in
$O(N\log N)$ ops using parallel computations.)

Let $ \zeta(z)=\sum_{k=0}^N\gm_kz^k. $

{\em Step 3.} Compute $\zeta_i(z)=[\wdt{\mba_{1i}}(z)\zeta(z)]^-$
for $j\not=i=1,2,\ldots,m$. (This step needs $ O(mN\log N)$ ops.)

{\em Step 4.} Use Algorithm 3.1 with the data
$\zeta_1(z),\zeta_2(z),\ldots,\zeta_{i-1}(z),\zeta_{i+1}(z),\ldots,\zeta_m(z)$
to construct the corresponding wavelet matrix. Denote this matrix
by $\mbA^\ddag(z)$ (in polyphase representation). Then, if we
transpose $\mbA^\ddag(z)$ and move its last row in the place of
$i$th row and its last column in the place of $i$th column, we get
the desired $\mbA(z)\in\calW_0(m,N,N)$.

All in all, the number of operations in Algorithms 3.1 and 3.2 can
be estimated as $O(m^2N\log N)+O(m^3N)$.

\section{Numerical Simulations}

To compare the performance of the described algorithms, their
computer code was written in Mathematica-8. A PC with 2.40GHz
Intel Quad Core CPU and 2GB RAM was used for numerical
simulations.

The accuracy level of the wavelet matrix construction algorithms
(2.1 and 3.1) is measured by that of relation (4), in which we
substitute the computation results, whereas the accuracy level of
wavelet matrix completion algorithms (2.2 and 3.2)  is naturally
measured by the difference between the initial data and the first
row of the computed matrix.

The accuracy levels  determined in this way are essentially the
same for both methods and are quite close to precisions in which
Mathematica-8 carries out calculations. (As it is known
Mathematica-8 provides an opportunity to make this precision
arbitrarily large.) Therefore, all experiments for comparison of
the speeds of different algorithms were run in the standard double
precision. The results of these simulations are presented in the
Tables below. As it was expected (since $O(m^2N^2)$ should be
larger than $O(m^2N\log N)+O(m^3N)$ for $m\ll N$), the new
algorithms are faster than the old ones, and the difference
between their performance times is becoming more and more evident
as $m$ and $N$ grow, keeping $N$ sufficiently larger than $m$ as
it should be for wavelet matrices applicable in practice.

\begin{center}
\small{Table I\\Results of Computer Simulations of Wavelet Matrix
Construction Algorithms}
\end{center}
{\fontsize{8}{8pt}\selectfont
\begin{center}
\begin{tabular}{|c|c|c|c|c|c|c|c|c|c|c|c|c|c|c|}
\hline &&&&&&&&&&&&&\\[-2mm]
Rank $m$&10&10&10&10&20&20&20&30&30&30&50&50&50
\\\hline &&&&&&&&&&&&&\\[-2mm]
Order $N$&50&150&300&400&100&250&300&100&150&200&100&200&300\\
\hline &&&&&&&&&&&&&\\[-2mm]
Time (New Alg.)&0.34&2.80&10.60&18.86&2.58&14.93&20.12&4.52&9.14&15.75&8.65&27.93&57.54\\
\hline &&&&&&&&&&&&&\\[-2mm]
Time (Old Alg.)&0.79&5.45&16.09&19.44&10.84&57.44&62.52&23.85&49.98&85.93&67.69&234.48&401.16\\
\hline
\end{tabular}
\end{center}}

\vskip+0.5cm

\begin{center}
\small{Table II\\Results of Computer Simulations of Wavelet Matrix
Completion Algorithms}
\end{center}
{\fontsize{8}{8pt}\selectfont
\begin{center}
\begin{tabular}{|c|c|c|c|c|c|c|c|c|c|c|c|c|c|c|}
\hline &&&&&&&&&&&&&\\[-2mm]
Rank $m$&10&10&10&10&20&20&20&30&30&30&50&50&50
\\\hline &&&&&&&&&&&&&\\[-2mm]
Order $N$&50&150&300&400&100&250&300&100&150&200&100&200&300\\
\hline &&&&&&&&&&&&&\\[-2mm]
Time (New Alg.)&0.39&2.83&10.96&19.76&2.85&15.76&22.77&4.58&9.46&17.13&9.23&28.01&58.48\\
\hline &&&&&&&&&&&&&\\[-2mm]
Time (Old Alg.)&0.85&6.08&18.24&22.01&11.96&63.13&70.54&24.96&52.91&88.09&74.42&333.71&473.61\\
\hline
\end{tabular}
\end{center}}

\section{Conclusion}

In this paper, we describe in detail two new algorithms of wavelet
matrix construction and completion, introduced in \cite{EL}. The
results of numerical simulations are presented, which prove the
advantage of the new algorithms over the existing algorithms in
performance speed.

Another advantage of the new method, which should be mentioned
here and might be used in the future, is that the algorithms based
on this method can be divided into $m$ parallel tasks which will
make them even more faster. The old method is heavily recurrent
and misses any opportunity to be parallelized.

 \vskip+0.2cm

\ Authors' Address: \vskip+0.2cm

\noindent  Faculty of Exact and Natural Sciences

\noindent I. Javakhishvili State University

\noindent 2, University Street, Tbilisi 0143, Georgia

\noindent E-mails: {\em \{nika.salia; alexander.gamkrelidze;
lasha.ephremidze\}@tsu.ge}

\end{document}